\documentclass[12pt,a4paper]{amsart}    
\usepackage[all]{xy}
\newtheorem{theorem}{Theorem}
\newtheorem{lemma}[theorem]{Lemma}
\newtheorem{prop}[theorem]{Proposition}
\newtheorem{cor}[theorem]{Corollary}
\theoremstyle{remark}
\newtheorem{example}[theorem]{Example}
\newtheorem{remark}[theorem]{Remark}

\newtheorem*{problem*}{Problem}
\newtheorem*{remark*}{Remark}
\newtheorem*{convention*}{Convention}
\newtheorem*{notation*}{Notation}
\newtheorem*{examples*}{Examples}
\newtheorem*{example*}{Example}
\newtheorem*{warning*}{Warning}
\def\N{{\mathbb N}}

\def\R{{\mathbb R}}
\def\T{{\mathbb T}}
\def\C{{\mathbb C}}
\def\Z{{\mathbb Z}}

\def\H{{\mathcal H}}
\def\K{{\mathcal K}}
\def\K{{\mathcal K}}
\def\B{{\mathcal B}}

\begin{document}
\title[Generalized Multiresolution Analyses]{Generalized multiresolution analyses\\ with given multiplicity functions}
\author[L. W. Baggett]{Lawrence~W.~Baggett}
\address{Lawrence Baggett, Department of Mathematics, University of Colorado, Boulder, Colorado 80309, USA}
\email{baggett@euclid.colorado.edu}
\author[N. S. Larsen]{Nadia~S.~Larsen}
\address{Nadia S. Larsen, Department of Mathematics, University of Oslo, Blindern, NO-0316 Oslo, Norway}
\email{nadiasl@math.uio.no}
\author[K. D. Merrill]{Kathy~D.~Merrill}
\address{Kathy Merrill, Department of Mathematics, Colorado College, Colorado Springs, Colorado, 80903, USA}
\email{kmerrill@coloradocollege.edu}
\author[J. A. Packer]{Judith~A.~Packer}
\address{Judith Packer, Department of Mathematics, University of Colorado, Boulder, Colorado 80309, USA}\email{packer@euclid.colorado.edu}

\author[I. Raeburn]{Iain~Raeburn}
\address{Iain Raeburn, School  of Mathematics and Applied Statistics, University of Wollongong, NSW 2522, Australia}
\email{raeburn@uow.edu.au}
\begin{abstract}
Generalized multiresolution analyses are increasing sequences of subspaces of a Hilbert space $\H$ that fail to be multiresolution analyses in the sense of wavelet theory because the core subspace does not have an orthonormal basis generated by a fixed scaling function.
Previous authors have studied a multiplicity function $m$ which, loosely speaking, measures the failure of the GMRA to be an MRA.
When the Hilbert space $\H$ is $L^2(\mathbb R^n)$, the possible multiplicity functions have been characterized by Baggett and Merrill.
Here we start with a function $m$ satisfying  a consistency condition which is known to be necessary, and build a GMRA in an abstract Hilbert space with multiplicity function $m$.
\end{abstract}

\thanks{Nadia Larsen was supported by the Research Council of Norway. Judith Packer was supported by the National Science Foundation. Iain Raeburn was supported by the Australian Research Council, through the ARC Centre for Complex Dynamic Systems and Control, and by the Centre de Recerca Matem{\` a}tica at the Universitat Aut{\` o}noma de Barcelona.}

\maketitle
\section*{Introduction}
A generalized multiresolution analysis (GMRA) for a Hilbert space $\H$ consists of an increasing sequence of closed subspaces $V_n$ such that the complements $W_n:=V_{n+1}\ominus V_n$ give a direct-sum decomposition $\H=\bigoplus_{n\in \Z}W_n$ in which the $W_n$ for $n\geq 0$ are invariant under a representation $\pi$ of an abelian group $\Gamma$  on $\H$, and in which $W_{n+1}$ is the dilation of $W_n$.
The representation theory of abelian groups associates to the representations $\pi|_{V_0}$ and $\pi|_{W_0}$ integer-valued multiplicity functions $m$ and $\widetilde m$ on the dual group $\widehat\Gamma$.
In this paper, we consider the question of which functions $m$ and $\widetilde m$ can arise as multiplicity functions of GMRAs.

Previous work on this question has focused on the case $\H=L^2(\mathbb R^n)$, with the group $\mathbb Z^n$ acting by translation and the dilation implemented by an integer matrix $A$ whose eigenvalues $\lambda$ satisfy $|\lambda|>1$.
In this case,  Baggett and Merrill showed that $m$ is associated to a GMRA if and only if $m$ satisfies a consistency condition, described in detail below, and a technical condition on the translates of the support of $m$ \cite{bm}, which was discovered independently by  Bownik, Rzeszotnik and Speegle \cite{BRS} in their characterization of the dimension function of a wavelet.
When a GMRA in $L^2(\mathbb R^n)$ has an associated (multi-) wavelet, one or more functions $\psi_k$ such that the translates $\pi_n(\psi_k)$ form an orthonormal basis for $W_0$ (so that $\widetilde m$ is constant), the  characterizations in \cite{bm} and \cite{BRS} coincide.

Here we show that the second technical condition from \cite{bm} and \cite{BRS} is particular to $L^2(\mathbb R^n)$: provided one is willing to consider GMRAs in abstract Hilbert spaces, there are surprisingly few restrictions on $m$ and $\widetilde m$ apart from the consistency condition of \cite{bm}.
Our new results include a general construction of filters for multiplicity functions (Proposition~\ref{constructH}), and a criterion for the purity of an associated isometry which improves a key result in \cite{bjmp} (Theorem~\ref{kathy}).

We begin in \S\ref{secmult} by discussing GMRAs and multiplicity functions, and stating our main theorem.
We work in an abstract Hilbert space, with a countable abelian group $\Gamma$ of translations and a dilation operator which is compatible with an endomorphism $\alpha$ of $\Gamma$.
In \S\ref{secfromLR}, we revisit the direct-limit construction from \cite{gpots} to see what extra input we need to ensure that the direct limit carries the necessary translation group and dilation operator (Theorem~\ref{absnonsense}).
Then in \S\ref{secproofmain} we prove our main theorem.
We first show that our multiplicity function $m$ admits a low-pass filter, which is a matrix $H$ of functions on $\widehat\Gamma$ satisfying relations, introduced in \cite{bjmp}, which generalize those of quadrature mirror filters.
From $H$ we build an isometry $S_H$ on a Hilbert space $\K$, following an idea which goes back at least to \cite{bj}, and Theorem~\ref{kathy} says that when the filter is low-pass, $S_H$ is a pure isometry.
Then, when we apply the construction of Theorem~\ref{absnonsense} to this isometry, we obtain a direct-limit Hilbert space which has the required GMRA.

Since we think Theorem~\ref{kathy} and its proof are likely to be of independent interest, we have made them the focus of a separate section.
 Our proof follows the general strategy suggested in \cite[Lemma 3.3]{bjmp}, but here we have been able to replace some of the grittier estimates with exact calculations, and those which remain are much sharper.
The crux of the argument is the almost everywhere pointwise convergence of a sequence of averages, which we achieve by applying the reverse martingale convergence theorem. In the final section we discuss some examples which show that our results have broader scope than those of \cite{bm} and \cite{BRS}.

\subsection*{Notation and standing assumptions} Throughout this paper, $\Gamma$ is a countable abelian group with compact dual $\widehat\Gamma$, and $\lambda$ denotes normalised Haar measure on $\widehat\Gamma$.
We fix an injective endomorphism $\alpha$  of $\Gamma$ such that $\alpha(\Gamma)$ has finite index $N$ in $\Gamma,$ 
and we write $\alpha^*$ for the endomorphism of $\widehat\Gamma$ onto itself defined by $\alpha^*(\omega) = \omega\circ \alpha$, and note that $|\ker\alpha^*|=N$.
We assume that  $\bigcup_{n\geq 1}\ker\alpha^{*n}$ is dense in $\widehat\Gamma$ (or equivalently, that $\bigcap_{n\geq 1}\alpha^n(\Gamma)=\{0\}$).

All Hilbert spaces in the paper are separable.

\section{Multiplicity functions and the main theorem}\label{secmult}

Let  $\pi:\Gamma\to U(\H)$ be a unitary representation, and let $\delta$ be a unitary operator on ${\mathcal H}$ such that
\[\delta^{-1}\pi_\gamma \delta = \pi_{\alpha(\gamma)}\ \text{ for all $\gamma\in\Gamma$.}\]
As in \cite{bjmp}, a \emph{generalized multiresolution analysis} (or GMRA) relative to $\pi$ and $\delta$ is a sequence $\{V_n:n\in {\mathbb Z}\}$ of closed subspaces
of ${\mathcal H}$ with the following properties:

\smallskip
\begin{itemize}
\item[(a)] $V_n\subset V_{n+1}$ for all $n$,

\smallskip
\item[(b)] $V_{n+1} = \delta(V_n)$ for all $n$,

\smallskip
\item[(c)] $\bigcup_{n=0}^\infty V_n$ is dense in ${\mathcal H}$ and $\bigcap_{n=-\infty}^0 V_n = \{0\}$, and

\smallskip
\item[(d)] $V_0$ is invariant under $\pi$.
\end{itemize}

\smallskip

Property (d) of a GMRA implies that $\rho:=\pi|_{V_0}$
is a unitary representation of $\Gamma$, and
Stone's theorem on unitary representations of abelian groups together with the multiplicity theory for projection valued measures \cite{george}
gives a Borel measure $\mu$ on $\widehat\Gamma$, unique Borel subsets $\sigma_1\supseteq \sigma_2\supseteq \ldots $ of $\widehat\Gamma,$
and a unitary operator $J:V_0\to \bigoplus_i L^2(\sigma_i,\mu)$
satisfying
$$[J(\rho_\gamma(v))]_i(\omega) = \omega(\gamma)[J(v)]_i(\omega)$$
for $v\in V_0$,  $\gamma\in\Gamma$ and $\mu$-almost all $\omega\in\widehat\Gamma$
(see \cite[Proposition~1]{bmm}).

When ${\mathcal H}=L^2({\mathbb R}^d)$,
$\Gamma$ is the lattice ${\mathbb Z}^d$,
$\alpha(k) = Ak$
and $\pi$ is the representation determined by translation, the measure $\mu$
is necessarily absolutely continuous with respect to the Haar measure on
the torus ${\mathbb T}^d \equiv \widehat{{\mathbb Z}^d}$ (see \cite[Propositions~2 and 3]{bmm}).
This absolute continuity does not necessarily hold in general, but here we are interested in the converse, and we assume that our measures $\mu$ are absolutely continuous with respect to the Haar measure $\lambda$ on $\widehat\Gamma$.

With the above conventions, the function $m=\sum \chi_{\sigma_i}$ is called the
\emph{multiplicity function} of the GMRA.
Properties (a), (b) and (d) in the definition of a GMRA
imply that the subspace $W_0=V_1\ominus V_0$ also is invariant under $\pi$,
and hence determines a unitary representation $\widetilde\rho$
of $\Gamma$ on $W_0$. As above, Stone's theorem gives a measure $\widetilde\mu$ on $\widehat\Gamma$,
subsets $\widetilde \sigma_1\supseteq \widetilde \sigma_2\supseteq\ldots$ of $\widehat\Gamma$,
and a unitary map $\widetilde J:W_0\to \bigoplus_k L^2(\widetilde \sigma_k,\widetilde\mu)$ such that
$$[\widetilde J(\widetilde\rho_\gamma(v))]_k(\omega) = \omega(\gamma)[\widetilde J(v)]_k(\omega)$$
for $v\in W_0$, $\gamma\in\Gamma$
and $\widetilde\mu$-almost all $\omega\in\widehat\Gamma$.
We write $\widetilde m$ for the corresponding \emph{complementary multiplicity
function} given by
\[\widetilde m(\omega) = \sum_k \chi_{\widetilde\sigma_k}(\omega).\]
We now prove that the multiplicity functions $m$ and $\widetilde m$ of the GMRA $\{V_n :n\in \Z\}$ satisfy the following consistency equation for $(\mu+\widetilde\mu)$-almost all $\omega\in \widehat\Gamma$:
\begin{equation}\label{consistency}
m(\omega) + \widetilde m(\omega) = \sum_{\alpha^*(\zeta)=\omega} m(\zeta).
\end{equation}

This will follow immediately from Lemma \ref{general-consistency} below, which gives the consistency equation under slightly more general assumptions reflecting the dependency between 
\eqref{consistency} and the properties (a), (b) and (d) alone in the definition of a GMRA. 

\begin{lemma}\label{general-consistency}
Let $\rho$ and $\tilde{\rho}$ be 
representations of $\Gamma$ on closed subspaces $V$ and $W$ of a Hilbert space 
$\mathcal{H}$ such that there is a unitary $\delta$ on $\mathcal{H}$ satisfying the following 
conditions:
\begin{enumerate}
\item[(i)] $\delta(V)=V\oplus W$, and
\item[(ii)] $\delta^{-1}(\rho\oplus \tilde{\rho})_\gamma\delta\vert_V =\rho_{\alpha(\gamma)}$ for all $\gamma \in 
\Gamma.$ 
\end{enumerate}
Let $m$, $\mu$ and $\widetilde m$, $\tilde{\mu}$ be the multiplicity functions and the associated 
Borel measures given by Stone's theorem for $\rho$ and respectively $\tilde{\rho}$. Then 
$
m(\omega) + \widetilde m(\omega) = \sum_{\alpha^*(\zeta)=\omega} m(\zeta) 
$ for $(\mu+\tilde{\mu})$-almost all $\omega \in \widehat{\Gamma}.$
\end{lemma}

\begin{proof}
We begin by recalling an additional consequence of Stone's theorem.  
Suppose that $\pi$ is a representation of the abelian group 
$\Gamma$ acting in a Hilbert space $\mathcal{V}$, and let a Borel measure $\nu$ and Borel subsets 
$\{\tau_i\}$ be as in the statement of Stone's theorem. 
Suppose $\{\tau'_l\}$ is another collection of (not necessarily nested) Borel
subsets of $\widehat\Gamma,$ and suppose $J^1$ is a unitary
operator from $\mathcal{V}$ onto $\bigoplus_l L^2(\tau'_l)$ satisfying
\[
[J^1(\pi_\gamma(f))](\omega) = \omega(\gamma) [J^1(f)](\omega) 
 \]
for all $f\in \mathcal{V},$ all $\gamma\in\Gamma,$ and $\nu$-almost all $\omega\in\widehat\Gamma.$
Then 
\[ \sum\chi_{\tau_i}(\omega) = \sum_l \chi_{\tau'_l}(\omega)\]
for $\nu$-almost all $\omega\in\widehat\Gamma$. (This is really part of the proof of Stone's theorem.
In fact, Stone's theorem is essentially the same as the canonical decomposition theorem 
for projection-valued measures; see \cite{george}.)

Let $\sigma_i$, $J$ and $\tilde{\sigma}_k$, $\tilde{J}$ be the nested Borel subsets 
and unitaries given by Stone's theorem applied to $\rho$ and respectively $\tilde{\rho}$. Let 
$m'$ be the multiplicity function associated to $\rho\oplus\tilde{\rho}$ on $V\oplus W$. 
Define a unitary $J^1$ from $V\oplus W$ to $\bigoplus_i L^2(\sigma_i) \oplus \bigoplus_k 
L^2(\widetilde\sigma_k)$ by
 \[
 [J^1(f\oplus  g)](\omega)= [J(f)](\omega) \oplus [\widetilde{J}(g)](\omega). \]
The additional consequence of Stone's theorem described above implies that for almost 
all $\omega$, we have 
 \[
m'(\omega)  = \sum_i\chi_{\sigma_i}(\omega) + \sum_k \chi_{\widetilde \sigma_k}(\omega) 
  = m(\omega) + \widetilde m(\omega).
 \]
 Thus to verify the equation claimed in the lemma, it suffices to prove that
\begin{equation}
m'(\omega) = \sum_{\alpha^*(\zeta)=\omega} m(\zeta).
\label{enough_for_consistency}
\end{equation}
Let $s$ be a Borel cross-section for the quotient map of $\widehat{\Gamma}$ onto 
$\widehat\Gamma/\ker(\alpha^*)$. For each $i$ and each $\eta$ in the kernel of $\alpha^*,$
define 
\[\tau_{i,\eta} = \{\omega\in \widehat{\Gamma}:s(\omega)\eta \in \sigma_i\}, \]
and using (i) define $J':V\oplus W\to \bigoplus_{i,\eta} L^2(\tau_{i,\eta})$ by
 \[[J'(f\oplus g)]_{i,\eta}(\omega) = [J(\delta^{-1}(f\oplus g))]_i(s(\omega)\eta). \]
For $f\in V$ and $g\in W$ we then have 
 \begin{align*}
 [J'((\rho\oplus \tilde{\rho})_\gamma(f\oplus g))]_{i,\eta}(\omega)
 & = [J(\delta^{-1}( (\rho\oplus \tilde{\rho})_\gamma(f\oplus g)))]_i(s(\omega)\eta) \cr
& = [J(\rho_{\alpha(\gamma)}(\delta^{-1}(f \oplus g )))]_i(s(\omega)\eta) \text{ by (ii) }\cr
& = [s(\omega)\eta](\alpha(\gamma)) [J(\delta^{-1}(f\oplus g))]_i(s(\omega)\eta) \cr
& = [\alpha^*(s(\omega)\eta)](\gamma) [J'(f\oplus g)]_{i,\eta}(\omega) \cr
& = \omega(\gamma) [J'(f\oplus g)]_{i,\eta}(\omega).
\end{align*}
Therefore, again by the additional consequence of Stone's theorem described above, the 
multiplicity function $m'$ is given by
\begin{align*}
m'(\omega) & = \sum_{i,\eta} \chi_{\tau_{i,\eta}}(\omega)  = \sum_{i,\eta} \chi_{\sigma_i}(s(\omega)\eta) \cr
& = \sum_\eta \sum_i \chi_{\sigma_i}(s(\omega)\eta)  = \sum_{\eta} m(s(\omega)\eta) \cr
& = \sum_{\alpha^*(\zeta)=\omega} m(\zeta),
 \end{align*}
as was sought in \eqref{enough_for_consistency}.
\end{proof}

\begin{remark}
We look at the above definitions in the familiar setting of $L^2(\R)$.   Suppose that $\pi$ is the representation of $\Gamma=\Z$ by translations on $L^2(\R)$ and $\delta$ is a dilation operator.
When $W_0$ is generated by a wavelet $\psi$, so that the translates $\{\pi_n\psi\}$ form an orthonormal basis for $W_0$, the representation $\pi|_{W_0}$ is equivalent to the representation by multiplication operators on $L^2(\T)$ and the complementary multiplicity function $\widetilde m$ is identically $1$.
When there is a scaling function $\phi$ such that $\{\pi_n\phi \}$ is an orthonormal basis for $V_0$, so that the GMRA is an MRA, we also have $m$ identically equal to $1$.
However, the Journ{\'e} wavelet provides an example of a wavelet such that the corresponding GMRA is not an MRA, and the multiplicity function $m$ is not constant.
The function $m$ for the Journ{\'e} wavelet is explicitly worked out in \cite{cou} (and see also Example~\ref{Journe} below).
\end{remark}

In this paper, we ask what functions can arise as multiplicity functions, and our main result is the following theorem.

\begin{theorem}
\label{achievem}
Suppose $c\in \N$ and $m:\widehat\Gamma \to \{0,1,\cdots,c\}$ is a Borel function which satisfies
the consistency inequality 
\begin{equation}\label{consistencyform}
m(\omega) \leq \sum_{\alpha^*(\zeta)=\omega} m(\zeta),
\end{equation}
and define $\widetilde m:\widehat\Gamma\to \{0,1,\cdots,c\}$ by
\begin{equation} \label{defmtilde}
\widetilde m(\omega) = \sum_{\alpha^*(\zeta)=\omega} m(\zeta) - m(\omega).
\end{equation}
Suppose that there  is a positive integer $a$ satisfying $m(\omega)-\tilde{m}(\omega)\leq a \leq m(\omega)$ for all $\omega$ near 1.  Then there is a GMRA which has $m$ and $\widetilde m$
as the associated multiplicity and complementary multiplicity functions.
\end{theorem} 

\begin{remark}
In many examples, the multiplicity function $m$ 
attains its maximum value $c$ throughout a neighborhood of 1,
and then $a=c$ satisfies the hypothesis
of the theorem.
\end{remark}
\section{Construction of a GMRA from a pure isometry}\label{secfromLR}

We extend some of the ideas from \cite{gpots} by proving the following theorem.

 \begin{theorem}\label{absnonsense}
 Suppose $S$ is an isometry on a Hilbert space ${\mathcal K}$, and let $(\K_\infty,U_n)$ be the direct limit of the direct system $(H_n,T_n)$ in which each Hilbert space $H_n={\mathcal K}$ and each $T_n=S$:

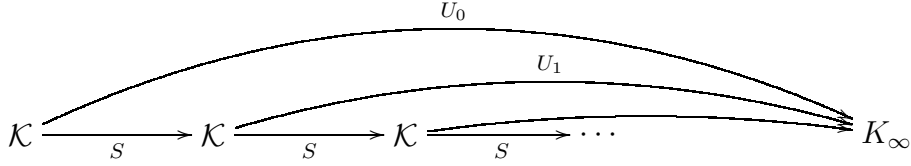
 \[\xygraph{{{\mathcal K}}="v0":[rr]{{\mathcal K}}="v1"_{S}:[rr]{{\mathcal K}}="v2"_{S}:[rr]{\cdots}="v3"_{S}:@{}[rrr]
{K_\infty}="v4"
 "v0":@/^{40pt}/^{U_0}"v4""v1":@/^{20pt}/"v4"^{U_1}"v2":@/^{7pt}/"v4"^{} }\]

\textnormal{(a)} There is a unitary operator $S_{\infty}$ on $\K_{\infty}$ such that
\[S_{\infty}(U_nh)=U_n(Sh)=U_{n-1}h\ \mbox{ for every $h\in H_n:={\mathcal K}$.}\]

 \textnormal{(b)} The subspaces $V_n:=U_n{\mathcal K}$ of $\K_{\infty}$ satisfy
 \begin{enumerate}
\item[(i)] $V_n\subset V_{n+1}$;
\smallskip
\item[(ii)] $\bigcup_{n=0}^\infty V_n$ is dense in $\K_{\infty}$;
\smallskip
\item[(iii)] $S_{\infty}$ is a unitary isomorphism of $V_{n+1}$ onto $V_n$.
\end{enumerate}

\textnormal{(c)} For $n<0$, define $V_n:=S_\infty^{|n|}(V_0)$.
Then $\bigcap_{n\in \Z}V_n=\{0\}$ if and only if $S$ is a pure isometry.

\smallskip
\textnormal{(d)}  If $\rho$ is a unitary representation of $\Gamma$
on ${\mathcal K}$ such that
 \begin{equation}\label{hyponrho} 
S\rho_\gamma = \rho_{\alpha(\gamma)} S\ \text{ for $\gamma\in\Gamma$},
\end{equation}
then there exists a unitary representation $\pi$ of $\Gamma$ on $\K_\infty$ such that all the subspaces $V_n,$ for $n\geq 0,$ are invariant under $\pi$ and
 \begin{equation}\label{conconpi}
S_\infty \pi_\gamma = \pi_{\alpha(\gamma)} S_\infty \ \text{ for $\gamma\in\Gamma$}.
 \end{equation}
 \end{theorem}

    \begin{proof}
The construction of $S_\infty$ is described on page~37 of \cite{gpots}, and it is proved there that $S_\infty$ is unitary. Since $U_{n}{\mathcal K}=U_{n+1}S{\mathcal K}$, we have $V_n\subset V_{n+1}$; that the union of the subspaces $V_n=U_n{\mathcal K}$ is dense is a standard property of the direct limit. The equation
\[ S_{\infty}(U_{n+1}h)=U_{n+1}(Sh)=U_nh \]
shows that $S_\infty$ is an isomorphism of $V_{n+1}$ onto $V_n$.

For (c), we notice that for $n<0$,
 \[V_n=S_{\infty}^{|n|}V_0=S_{\infty}^{|n|}U_0{\mathcal K}=U_0S^{|n|}{\mathcal K},\]
so
\begin{align*}
 \textstyle{\bigcap_{n\in \Z}V_n=\{0\}}&\Longleftrightarrow\textstyle{\bigcap_{k=1}^{\infty}V_{-k}=\{0\}}\\
&\Longleftrightarrow\textstyle{\bigcap_{k=1}^{\infty}U_0S^k{\mathcal K}=\{0\}}\\
 &\Longleftrightarrow\textstyle{\bigcap_{k=1}^{\infty}S^k{\mathcal K}=\{0\}}.
 \end{align*}
 Since $\bigcap_{k=1}^{\infty}S^k{\mathcal K}$ is the largest subspace of ${\mathcal K}$ on which $S$ is unitary, this proves (c).

The intertwining relation \eqref{hyponrho} implies that 
\[\xygraph{
 {\K}="v0":[rr]{\K}="v1"^{S}:[rr]{\K}="v2"^{S}:[rr]{\cdots}="v3"^{S}:@{}[r]{\cdots}="v_3":[rr]
 {\K_\infty}="v4"
 "v0":[dd]{\K}="w0"^{\rho_\gamma}
"w0":[rr]{\K}="w1"^{S}:[rr]{\K}="w2"^{S}:[rr]{\cdots}="w3"^{S}:@{}[r]{\cdots}="w_3":[rr]
 {\K_\infty,}="w4"
 "v1":"w1"^{\rho_{\alpha(\gamma)}}"v2":"w2"^{\rho_{\alpha^2(\gamma)}}"v4":@{-->}"w4"_{\pi_\gamma}
 }\]
is a commutative diagram of isometries, and hence the universal property of the direct limit gives the existence of a unique isometry $\pi_\gamma$ such that 
 $\pi_\gamma\circ U_n=U_n\circ \rho_{\alpha^n(\gamma)}$, which implies immediately that $V_n=U_n{\mathcal K}$ is invariant under $\pi_\gamma$.
 Uniqueness implies that $\pi_{-\gamma}$ is an inverse for $\pi_\gamma$ and that $\pi_{\gamma\tau}=\pi_\gamma\pi_\tau$, so $\pi$ is a unitary representation of $\Gamma$.
 Finally, we have
 \begin{align*}
 S_\infty\pi_\gamma U_n&=S_\infty U_n\rho_{\alpha^n(\gamma)}=U_nS\rho_{\alpha^n(\gamma)}\\
 &=U_n\rho_{\alpha^{n+1}(\gamma)}S=\pi_{\alpha(\gamma)}U_nS\\
 &=\pi_{\alpha(\gamma)}S_\infty U_n,
 \end{align*}
 which establishes \eqref{conconpi}.
  \end{proof}

   \begin{cor}
If $S$ is a pure isometry on $\K$, then the subspaces $V_n$ of $\K_\infty$ form a generalized multiresolution analysis with respect to $\pi:\Gamma\to U(\K_\infty)$ and $\delta:=S_\infty^{-1}$.
 \end{cor}

 \section{Proof of the main theorem}\label{secproofmain}
 
 Let $m:\widehat\Gamma\to \Z$ be a Borel function such that $0\leq m(\omega)\leq c$ for all $\omega$, and for $0\leq i\leq c$ write $\sigma_i:=\{\omega\in\widehat\Gamma:m(\omega)\geq i\}$.
  A \emph{filter relative to $m$ and $\alpha^*$} is a Borel function $H=[h_{i,j}]:\widehat\Gamma\to M_c(\C)$ such that $h_{i,j}$ vanishes outside $\sigma_j$ and 
  \begin{equation}\label{origfilter}
 \sum_{\alpha^*(\zeta)=\omega} H(\zeta)H^*(\zeta) =N\Sigma(\omega) \ \text{ for almost all $\omega\in \widehat\Gamma$},
 \end{equation}
where $\Sigma(\omega)$ is the diagonal matrix with entries $\chi_{\sigma_i}(\omega)$.
 Such a filter is \emph{low-pass of rank $a$} if $H$ is continuous near $1$ and $H(1)$ has block form
 \[
 H(1)=\begin{pmatrix}N^{1/2}1_a&0\\0&0\end{pmatrix}.
 \]

Crucial for our argument is that, when $m$ satisfies the hypotheses of Theorem~\ref{achievem}, there are always compatible low-pass filters.

\begin{prop}\label{constructH}
 Suppose that the positive integer $a$ satisfies
 \begin{equation}\label{condona}
 m(\omega)-\widetilde m(\omega)\leq a\leq m(\omega)\ \text{ for all $\omega$ near $1$ in $\widehat\Gamma$.
}
\end{equation}
Then there is a filter $H$ relative to $m$ and $\alpha^*$ which is low-pass of rank $a$.
 \end{prop}

 \begin{proof}
 We begin by writing the filter equations \eqref{origfilter} in the form
 \begin{equation}\label{filter}
 \sum_j \sum_{\zeta \in\ker \alpha^*} h_{i,j}(\omega\zeta) 
 \overline{h_{i',j}(\omega\zeta)}
 = N\delta_{i,i'}\chi_{\sigma_i}(\alpha^*(\omega)).
\end{equation}

    We choose a Borel cross-section $s$ for $\alpha^*$, and write $C=s(\widehat\Gamma)$; then every element in $\widehat\Gamma$ can be written in a unique way as $\omega\zeta$ for some $\omega\in C$ and $\zeta\in\ker\alpha^*$, and to build a filter it suffices to construct $c$ functions
 \[
 h_i(\omega)=\{h_{i,j}(\omega\zeta):1\leq j\leq m(\omega),\ \zeta\in \ker\alpha^*\}
 \]
from $C$ to $\C^{cN}$ such that the $h_{i,j}(\omega\zeta)$ vanish unless $\omega\zeta\in \sigma_j$ and \eqref{filter} holds for every $\omega\in C$.
 Equation~\eqref{filter} is equivalent to asking that the  vectors $(h_i(\omega))$ in $\C^{cN}$are orthogonal of norm $N^{1/2}\chi_{\sigma_i}(\alpha^*(\omega))$.

  Let $U$ be a neighborhood of $1$ such that \eqref{condona} holds for $\omega\in U$, and shrink $U$ to ensure that the sets $\{U\zeta:\zeta\in\ker\alpha^*\}$
are pairwise disjoint.

From the continuity of $\alpha^*$, there exist neighborhoods $V$
and $W$ of the identity, both contained in $U$, such that $\alpha^*$ maps
$W$ onto 
$V$, and since $U\cap U\zeta=\emptyset$ for $\zeta\not=1$, $\alpha^*$ is a homeomorphism of $W$ onto $V$.
We may suppose without loss of generality that $W\subset C$.    

 For $\omega\in W$ and $i\leq a$, we define $h_i(\omega)$ by
 \[
 h_{i, j}(\omega\zeta)= \begin{cases}
 \sqrt{N} &\text{if $i=j\leq a$ and $\zeta=1$}\\
 0 & \text{otherwise.
 }
 \end{cases}
 \]
 This will ensure that our filter is continuous at $1$ and is low-pass of rank $a$.
 For $i>m(\alpha^*(\omega))$, we must set $h_i(\omega)=0$ for all $\omega$.
 For $a<i\leq m(\alpha^*(\omega))$, the entries $h_{i,j}(\omega)$  must be $0$ for $\omega\in W$.   We also need to take $h_{i,j}(w\zeta)=0$ unless $w\zeta\in \sigma_j$, which is equivalent to $m(\omega\zeta)\geq j$.  Thus for each $\zeta$, there are $m(\omega\zeta)$ $j$s for which $h_{i,j}(w\zeta)$ can be non-zero, and hence $\sum_{\zeta\not= 1}m(\omega\zeta)$ potentially non-zero elements.

Since $\omega\in W$ implies $\alpha^*(\omega)\in V$, and since $V$ is contained in $U$, we have
 \[
m(\alpha^*(\omega))-a \leq \widetilde m(\alpha^*(\omega))= \sum_{\zeta \in \ker \alpha^*} m(\omega\zeta) - m(\omega) = \sum_{\zeta\in\ker \alpha^*,\;\zeta\neq 1} m(\omega\zeta).
 \]
Thus the number of components in $h_i(\omega)$ which can be non-zero is greater than or equal to the required number $m(\alpha^*(\omega))-a$ of orthogonal vectors $h_i(\omega)$, and it is possible to find such vectors.
 Since there are only finitely many possible sets of values of $m(\alpha^*(\omega))$ and $m(\omega\zeta)$, and we can use the same vectors for $h_i(\omega)$ when these values are all the same, we can find simple functions $h_i$ with the required properties.

Defining the vectors $h_i(\omega)$ for $\omega\in C\backslash W$ is easier, since now we just need to define $h_i(\omega)$ for $i\leq m(\alpha^*(\omega))$, and we have $\sum_{\zeta\in\ker\alpha^*}m(\omega\zeta)\geq m(\alpha^*(\omega))$ non-zero entries to play with.
 \end{proof}

Our main technical result shows that low-pass filters give rise to pure isometries.
In an attempt to clarify our overall strategy, we will postpone the proof of this result till the next section.

\begin{theorem}
  \label{kathy}
 Suppose that $m:\widehat\Gamma\to \{0,1,\cdots,c\}$ is Borel, and that $H$ is a filter relative to $m$ and $\alpha^*$. Let ${\mathcal K}$ be the Hilbert space defined by
\begin{equation}
 {\mathcal K} = \bigoplus_i L^2(\sigma_i),
 \end{equation}
and define an operator $S_H$ on $\mathcal{K}$ by
\begin{equation}
(S_Hf)(\omega) = H^t(\omega)f(\alpha^*(\omega)).
 \end{equation}
Then $S_H$ is an isometry on ${\mathcal K}$. If the filter $H$ is low-pass of some rank $a$ between $1$ and $c$, then $S_H$ is a pure isometry.
\end{theorem}

    We now have all the ingredients to prove our main theorem.

\begin{proof}[Proof of Theorem~\ref{achievem}]
Proposition~\ref{constructH} gives us a low-pass filter $H$ of rank $a$. Let $S_H$ be the pure isometry on ${\mathcal K} = \bigoplus_i L^2(\sigma_i)$
discussed in Theorem \ref{kathy}. Define a representation $\rho$ of $\Gamma$ on ${\mathcal K}$ by
\[\rho_\gamma(f)(\omega) = \omega(\gamma)f(\omega).
\]
Then
 \begin{align*}
 S_H(\rho_\gamma(f))(\omega)
  & = H^t(\omega) \rho_\gamma(f)(\alpha^*(\omega)) \cr
 & = H^t(\omega)\alpha^*(\omega)(\gamma) f(\alpha^*(\omega)) \cr
 & = H^t(\omega)\omega(\alpha(\gamma)) f(\alpha^*(\omega)) \cr
& = \rho_{\alpha(\gamma)}(S_Hf)(\omega).
 \end{align*}
Now applying Theorem~\ref{absnonsense} gives us a direct limit Hilbert space $(\K_\infty, U_n)$, a representation $\pi:\Gamma\to U(\K_\infty)$, and a dilation operator $\delta=S_\infty^{-1}$, such that $\{V_n\}:=\{U_n\K\}$
is a GMRA relative to $\pi$ and $\delta$.

 The canonical embedding $U_0$ is an isomorphism of ${\mathcal K} = \bigoplus_i L^2(\sigma_i)$ onto $V_0$ which intertwines $\rho$ and $\pi|_{V_0}$, so this 
 GMRA has multiplicity function $m$. It follows from equation~\eqref{consistency} that $\widetilde m$ must be the complementary multiplicity function.
\end{proof}

 This completes the proof of Theorem~\ref{achievem}, modulo our obligation to provide a proof of Theorem~\ref{kathy}.

\section{Low-pass filters and pure isometries}

In this section, we pay our debts by proving Theorem~\ref{kathy}.

 Since $H^t(\omega)_{ij}=h_{j,i}(\omega)$ has support in $\sigma_i$, $S_Hf$ belongs to $\K$.
A computation shows that the adjoint of $S_H$ is given by
 \[
(S_H^*f)(\omega) = \frac1N \sum_{\alpha^*(\zeta)=\omega} \overline{H(\zeta)} f(\zeta),
 \]
and it is then easy to check that $S_H^*S_Hf=f$ when $f$ belongs to $\K$.
 So $S_H$ is an isometry.

From here we assume that $H$ is low-pass, and aim to prove that $S_H$ is a pure isometry, or in other words that $\bigcap_{n=0}^\infty S_H^n\K=\{0\}$.
We assume that this is not true, and look for a contradiction.
 Since every non-zero Hilbert space contains a unit vector, we can find a unit vector $f$ in $\bigcap_{n=0}^\infty S_H^n\K$.
To arrive at our contradiction, we consider the sequence $f_n:=S_H^{*n}f$; since $S_H$ is unitary on $\bigcap_{n=0}^\infty S_H^n\K$ with inverse $S_H^*$, $\{f_n\}$ is a sequence of unit vectors in $\bigcap_{n=0}^\infty S_H^n\K$.

We will need to deal with the powers of $S_H$ and $S_H^*$, and simple induction arguments yield the following explicit formulas:
\begin{align*}
(S_H^nf)(\omega)&= \Big(\prod_{k=0}^{n-1} H^t({\alpha^*}^k(\omega))\Big) f({\alpha^*}^n(\omega)),\ \text{ and}\label{impSn}\\
(S_H^{*n}f)(\omega)&= \frac1{N^n} \sum_{{\alpha^*}^n(\zeta)=\omega}
\Big(\prod_{k=n-1}^0 \overline{H({\alpha^*}^k(\zeta))}\Big)f(\zeta).
\end{align*}

We view elements $g\in\K$ as functions from $\widehat\Gamma$ to $\C^c$ whose $i$th coordinate $g_i$ has support in $\sigma_i$, and write $\|g(\zeta)\|$ for the norm of the vector $g(\zeta)\in\C^c$.
Then for each $g\in \K$, the function $\zeta\mapsto \|g(\zeta)\|^2$ is integrable on $\widehat\Gamma$.
 We want to identify the integrable functions associated to our sequence $f_n=S_H^{*n}f$.
 A crucial step in the calculation is the following extension of the filter identity~\eqref{origfilter}.
 In the following formula~\eqref{superfilter} it is crucial that the products are interpreted in the correct order: the middle terms, for example, are the ones for which $k=0$ and~$l=0$.

\begin{lemma}
For every $n\geq 1$, we have
 \begin{equation}\label{superfilter}
 \sum_{\alpha^{*n}(\zeta)=\omega} \Big(\prod_{k=n-1}^0H(\alpha^{*k}(\zeta))\Big)\Big(\prod_{l=0}^{n-1}H^*(\alpha^{*l}(\zeta))\Big)= N^n\Sigma(\omega) \ \text{ for almost all $\omega$}.
 \end{equation}
\end{lemma}

 \begin{proof}
For $n=1$ we recover the usual filter identity \eqref{origfilter}.
Suppose \eqref{superfilter} is true for $n\geq 1$.
Then
\begin{align*}
 \sum_{\alpha^{*(n+1)}(\zeta)=\omega} \Big(\prod_{k=n}^0&H(\alpha^{*k}(\zeta))\Big) \Big(\prod_{l=0}^{n}H^*(\alpha^{*l}(\zeta))\Big)\\
&=\sum_{\alpha^{*n}(\eta)=\omega}\;\sum_{\alpha^*(\zeta)=\eta} \Big(\prod_{k=n}^1H(\alpha^{*k}(\zeta))\Big)H(\zeta)H^*(\zeta)\Big(\prod_{l=1}^{n}H^*(\alpha^{*l}(\zeta))\Big)\\
 &=\sum_{\alpha^{*n}(\eta)=\omega} \Big(\prod_{k=n-1}^0H(\alpha^{*k}(\eta))\Big)\Big(\sum_{\alpha^*(\zeta)=\eta}H(\zeta)H^*(\zeta)\Big)\Big(\prod_{l=0}^{n-1}H^*(\alpha^{*l}(\eta))\Big)\\
 &=\sum_{\alpha^{*n}(\eta)=\omega} \Big(\prod_{k=n-1}^0H(\alpha^{*k}(\eta))\Big)N\Sigma(\eta)\Big(\prod_{l=0}^{n-1}H^*(\alpha^{*l}(\eta))\Big).
 \end{align*}
 Now notice that for each $i$, both the $(i,j)$ entry $h_{i,j}(\eta)$ in $H(\eta)$ and the $(j,i)$ entry in $H^*(\eta)$ vanish unless $\eta\in \sigma_j$, in which case the $(j,j)$ entry in $\Sigma(\eta)$ is $1$. So the $\Sigma(\eta)$ in the middle has no effect, and we deduce from the inductive hypothesis that the last expression reduces to $N(N^n\Sigma(\omega))=N^{n+1}\Sigma(\omega)$.
 \end{proof}

 \begin{lemma}\label{weeglimmer}
 For almost all $\omega$ we have
 \begin{equation}\label{exactvofest}
 \|f_n(\omega)\|^2=\frac1{N^n} \sum_{{\alpha^*}^n(\zeta) = \omega}
    \|f(\zeta)\|^2.
\end{equation}
 \end{lemma} 

\begin{proof} 
We write $(v\,|\,w)$ for the usual inner product on $\C^c$.
Then since $f\in S_H^n\K=S_H^nS_H^{*n}\K$, we have
\begin{align}\label{calcnormfn}
 \frac1{N^n} \sum_{{\alpha^*}^n(\zeta) = \omega}
 &\|f(\zeta)\|^2
 =\frac1{N^n} \sum_{{\alpha^*}^n(\zeta) = \omega}(f(\zeta)\,|\,f(\zeta))\\
&=\frac1{N^n} \sum_{{\alpha^*}^n(\zeta) = \omega}(S_H^nS_H^{*n}f(\zeta)\,|\,S_H^nS_H^{*n}f(\zeta))\notag\\
&=\frac1{N^n} \sum_{{\alpha^*}^n(\zeta) = \omega}(S_H^nf_n(\zeta)\,|\,S_H^nf_n(\zeta))\notag\\
 &=\frac1{N^n} \sum_{{\alpha^*}^n(\zeta) = \omega}\Big(\prod_{l=0}^{n-1}H^t(\alpha^{*l}(\zeta))f_n(\alpha^{*n}\zeta)\,\Big|\,\prod_{k=0}^{n-1}H^t(\alpha^{*k}(\zeta))f_n(\alpha^{*n}\zeta)\Big)\notag\\
 &=\frac1{N^n} \sum_{{\alpha^*}^n(\zeta) = \omega}\bigg(\Big(\prod_{k=n-1}^0\overline{H(\alpha^{*k}(\zeta))}\Big)\Big(\prod_{l=0}^{n-1}H^t(\alpha^{*l}(\zeta))\Big)f_n(\omega)\,\Big|\,f_n(\omega)\bigg)\notag\\
    &=(\Sigma(\omega)f_n(\omega)\,|\,f_n(\omega)),\notag
 \end{align}
 where at the last step we used the conjugate of \eqref{superfilter}.
Now we deduce from the original filter equation \eqref{origfilter} that
 \[
 h_{i,j}(\zeta)\not=0\Longrightarrow \alpha^*(\zeta)\in \sigma_i.
 \]
Then the $i$th entry $[f_n(\omega)]_i$ satisfies
 \[
 [f_n(\omega)]_i=\frac{1}{N^n}\sum_{\alpha^{*n}(\zeta)=\omega}\sum_{j=i}^c \overline{h_{i,j}(\alpha^{*(n-1)}(\zeta))}\Big[\Big(\prod_{k=n-2}^0\overline{H(\alpha^{*k}(\zeta))}\Big)f(\zeta)\Big]_j,
  \]
and hence vanishes unless $\alpha^*(\alpha^{*(n-1)}(\zeta))=\omega$ is in $\sigma_i$.
Thus $\Sigma(\omega)f_n(\omega)=f_n(\omega)$, and the calculation \eqref{calcnormfn} gives the result.
 \end{proof}

We can rewrite the formula \eqref{exactvofest} as
 \begin{equation}\label{exactvofest2}
 \|f_n(\alpha^{*n}(\omega))\|^2=\frac1{N^n} \sum_{\eta \,\in\,\ker\alpha^{*n}}
  \|f(\omega\eta)\|^2,
 \end{equation}
and now we claim that the right-hand side $X_n(\omega)$ of \eqref{exactvofest2} is the expectation $E(\|f\|^2\,|\,\B_n)$ of $f$ with respect to the subalgebra $\B_n:=(\alpha^*)^{-n}(\B)$ of the Borel $\sigma$-algebra $\B$.
To see this, we note that $\B_n$ is the $\sigma$-algebra of Borel sets which are invariant under the action of $\ker \alpha^{*n}$, so that $X_n$ is certainly $\B_n$-measurable, and for $B\in \B_n$, we have
\begin{align*}
 \int_B X_n(\omega)\,d\omega&=\frac{1}{N^n}\sum_{\zeta\,\in\,\ker\alpha^{*n}}\int_{B\zeta^{-1}} \|f(\omega)\|^2\,d\omega\\&=\frac{1}{N^n}\sum_{\zeta\,\in\,\ker\alpha^{*n}}\int_{B} \|f(\omega)\|^2\,d\omega\\
 &=\int_B \|f(\omega)\|^2\,d\omega;
 \end{align*}
 in other words, $X_n$ has the properties which characterise $E(\|f\|^2\,|\,\B_n)$ (see the observation at the top of page~18 of \cite{pet}), and hence $X_n=E(\|f\|^2\,|\,\B_n)$.
Further, if $B\in \B_{n+1}\subset\B_n$, then we have
 \[
 \int_B X_{n+1}(\omega)\,d\omega=\int_B \|f(\omega)\|^2\,d\omega=\int_B X_{n}(\omega)\,d\omega,
\]
and hence $X_{n+1}=E(X_n\,|\,\B_{n+1})$.
Thus the family $\{X_n\}$ satisfies the hypotheses of the reverse martingale convergence theorem (as in \cite[Theorem~10.
6.1]{dud}, for example), and we can deduce from that theorem that $X_n$ converges almost everywhere to the expectation $E(\|f\|^2\,|\,\B_\infty)$ associated to $\B_\infty:=\bigcap_{n\geq 1}\B_n$.

To identify $E(\|f\|^2\,|\,\B_\infty)$, we need the following standard lemma.

\begin{lemma}
 If $B\in\B_\infty$, then $\lambda(B)$ is either $0$ or $1$.
 \end{lemma}

 \begin{proof}
 Notice that $B$ is invariant under multiplication by elements of $\ker\alpha^{*n}$ for every $n\geq 1$.
Suppose $\gamma\in \Gamma\backslash\{0\}$.
 Since $\bigcup_{n\geq 1}\ker\alpha^{*n}$ is dense in $\widehat\Gamma$, two characters of $\widehat\Gamma$ which agree on $\bigcup_{n\geq 1}\ker\alpha^{*n}$ must agree on all of $\widehat\Gamma$.
 Thus there exist $n$ and $\zeta\in\ker\alpha^{*n}$ such that $\zeta(\gamma)\not=1$.
Then the Fourier coefficients of the characteristic function $\chi_B$ satisfy
\begin{align*}
 \widehat \chi_B(\gamma)&=\int_{\widehat\Gamma}\chi_B(\omega)\omega(\gamma)\,d\omega=
  \int_{\widehat\Gamma}\chi_B(\zeta\omega)(\zeta\omega)(\gamma)\,d\omega\\
 &=\zeta(\gamma)\int_{\widehat\Gamma}\chi_{\zeta^{-1}B}(\omega)\omega(\gamma)\,d\omega=\zeta(\gamma)\int_{\widehat\Gamma}\chi_{B}(\omega)\omega(\gamma)\,d\omega\\
&=\zeta(\gamma)\widehat \chi_B(\gamma),
 \end{align*}
 and hence $\widehat \chi_B(\gamma)=0$.
Thus $\widehat \chi_B(\gamma)=0$ for every non-zero $\gamma$, and $\chi_B$ is either $0$ or $1$ in $L^1(\widehat\Gamma)$, which implies the result.
\end{proof}

 So $\B_\infty=\{B\in\B:\lambda(B)=0\text{ or }1\}$, and the expectation $E(\|f\|^2\,|\,\B_\infty)$ is the constant function $\int \|f(\omega)\|^2\,d\omega$.
Since our $f$ is a unit vector, we have now proved the following Proposition.

  \begin{prop} For almost all $\omega\in \widehat\Gamma$, we have
\[ \|f_n(\alpha^{*n})(\omega)\|^2=\frac1{N^n} \sum_{\eta \,\in\,\ker\alpha^{*n}}
\|f(\omega\eta)\|^2\to 1\ \text{ as $n\to\infty$}. \]
\end{prop} 

We are now ready to get the contradiction which will prove that $S_H$ is a pure isometry.
We fix $\delta>0$.
 We view $H$ as a block matrix $H=(H_{i,j})$ for the decomposition $\C^c=\C^a\oplus \C^{c-a}$,
 and choose a neighborhood $V$ of the identity such that
 for each $\omega\in V$ we have
 \[
 \|H_{1,1}(\omega)-\sqrt N1_a\| < \delta\ \text{ and }\ \|H_{i,j}(\omega)\| < \delta \text{ for $(i,j)\not=(1,1)$.} \]
 Next, we choose a neighborhood $W$ of the identity such that $W$,
${\alpha^*}(W)$ and ${\alpha^*}^2(W)$ are all contained in $V$.

  By Egorov's theorem, there exists a set $E$ whose complement
has measure less than $\lambda(W)/4$ and an integer $M$ such that,
for all $n\geq M$ and all $\omega\in E$,
 \[1-\delta < \|f_n(\alpha^{*n}(\omega))\| < 1+\delta. \]
 Since $\alpha^*$ is measure-preserving in the sense that $\lambda((\alpha^*)^{-1}(E))=\lambda(E)$, the set 
  \[
 A:=W\cap (\alpha^*)^{-M}(E)\cap (\alpha^*)^{-(M+1)}(E)\cap (\alpha^*)^{-(M+2)}(E)
 \]
has positive measure.
 It then follows 
that $\alpha^{*M}(\omega),$  $\alpha^{*(M+1)}(\omega),$
  and $\alpha^{*(M+2)}(\omega)$
 all belong to $V\cap E$ for every $\omega\in A$.
 \par

We now fix $\omega\in A$, write $v=([v]_1,[v]_2)$ for the block decomposition of $v\in\C^c$, and make lower and upper estimates for $\|[f_{M+1}(\alpha^{*(M+1)}(\omega))]_1\|$. For the lower estimate, we observe that 
\begin{align*}
 \|[f_{M+1}(\alpha^{*(M+1)}(\omega))]_2\|
 &\leq\sum_{j=1}^2\big\|H^t_{2,j}(\alpha^{*(M+1)}(\omega)) [f_{M+2}(\alpha^{*(M+2)}(\omega))]_j\big\|\\
& \leq 2\delta (1+\delta),
\end{align*}
and deduce that
 \begin{align*}
 \|[f_{M+1}(\alpha^{*(M+1)}(\omega))]_1\| 
 & \geq \big\|f_{M+1}(\alpha^{*(M+1)}(\omega))\big\|- \big\|[f_{M+1}(\alpha^{*(M+1)}(\omega))]_2\big\|\\
& \geq 1-\delta -2\delta(1+\delta).
 \end{align*}
For the upper estimate, we write
\begin{equation}\label{2ndgo}
[f_{M}(\alpha^{*M}(\omega))]_1 =\sum_{j=1}^2H_{1,j}(\alpha^{*M}(\omega)) [f_{M+1}(\alpha^{*(M+1)}(\omega))]_j,
\end{equation}
 rewrite the first summand on the right as
  \[
 N^{1/2} [f_{M+1}(\alpha^{*(M+1)}(\omega))]_1
 +  (H_{1,1}(\alpha^{*M}(\omega))-N^{1/2}1_a) [f_{M+1}(\alpha^{*(M+1)}(\omega))]_1,
 \]
 and turn \eqref{2ndgo} round to get the estimate
 \[
 N^{1/2}\big\|[f_{M+1}(\alpha^{*(M+1)}(\omega))]_1\big\|
 \leq 1+\delta+2\delta(1+\delta)=(1+\delta)(1+2\delta).
 \]
Combining the upper and lower estimates shows that for every $\delta>0$ we must have
\[
1-\delta -2\delta(1+\delta) \leq
N^{-1/2}(1+\delta)(1+2\delta),
\]
which we can see is impossible by letting $\delta\to 0$.

This completes the proof of Theorem~\ref{kathy}.

\section{Examples}   

  We give several examples showing how low-pass filters can occur for different values of $a$ in Proposition \ref{constructH} and Theorem \ref{kathy} even for the same multiplicity function $m$.  For simplicity of notation, we identify the multiplicative group $\mathbb T$ with the additive set $\mathbb R/\mathbb Z,$ where we choose the coset representatives of $\mathbb R/\mathbb Z$ in $\mathbb R$ to be $[-1/2,1/2)$. This agrees with the notation in \cite{bjmp}.

\begin{example}\label{Journe}
Consider the multiplicity function for dilation by $2$ in $\mathbb R$ corresponding to the Journ\'e wavelet, previously studied in \cite{cou}, \cite{bcm} and \cite{bjmp}.  The multiplicity function for this minimally supported frequency wavelet is given by 
\[
m(x)=\begin{cases}2&\text{if $x\in [-\frac{1}{7},\frac{1}{7})$}\\
1&\text{if  $x\in \pm [\frac{1}{7},\frac{2}{7})\cup \pm[\frac{3}{7},\frac{1}{2})$}\\
0&\text{otherwise.}
\end{cases}
\]
Filters which give rise to the Journ\'e wavelet and satisfy the low-pass condition of rank $a=1$ were constructed in \cite{cou}.  However, for this $m$ the number $a=2$ also satisfies the conditions of Proposition \ref{constructH}, and we can also find filters which satisfy the low-pass condition of rank $a=2$.  Indeed, 
\[
h_{1,1}=\sqrt{2}\chi_{[-\frac{2}{7},-\frac{1}{4})\cup [-\frac{1}{7},\frac{1}{7})\cup [\frac{1}{4}, \frac{2}{7})},\quad  h_{1,2}=h_{2,1}= 0, \quad \text{and}\quad h_{2,2}=\sqrt{2}\chi_{[-\frac{1}{14},\frac{1}{14})}
\] 
have the required properties. If we consider the $2\times 2$ matrix $H=[h_{i,j}]$ as in Theorem~\ref{kathy} and view the infinite product $\prod_{j=1}^{\infty}[2^{-1/2}H(2^{-j}x)]$ as an element of $M(2,L^{\infty}(\mathbb R))$, we obtain a diagonal matrix with diagonal entries \[
\hat{\phi}_1:=\chi_{[-4/7,-1/2)\cup [-2/7,2/7)\cup [1/2,4/7)}\quad\text{ and}\quad\hat{\phi}_2:=\chi_{[-\frac{1}{7},\frac{1}{7})}.
\]
Then the shift invariant subspace of $L^2(\R)$ generated by the inverse Fourier transforms  $\phi_1$ and $\phi_2$ has multiplicity function 
\[
m'=\chi_{[-\frac{1}{2},-\frac{3}{7})\cup [-\frac{2}{7},\frac{2}{7})\cup [\frac{3}{7},\frac{1}{2})},
\]
which is a degenerate version of the original multiplicity function $m$. Thus one does not have an analogue of Theorem 3.4 of \cite{bcm} in this situation, and one must rely on the direct limit Hilbert space instead of $L^2(\mathbb R)$ to construct our GMRA.
\end{example}

\begin{example} The Journ\'e wavelet set actually corresponds to the case $n=1$ \cite[\S4, Example 1]{bmm}).  If we consider the case $n=2$  in this family of examples, we obtain the wavelet set $W = \pm [\frac{4}{15},\frac{1}{2}) \cup \pm[4, \frac{64}{15})$ and the multiplicity function given by \[
m(x)=\begin{cases}3&\text{if $x\in [-\frac{1}{15},\frac{1}{15})$}\\
2&\text{ if $x\in \pm [\frac{1}{15},\frac{2}{15})$}\\
1&\text{if  $x\in \pm [\frac{2}{15},\frac{4}{15})\cup \pm [\frac{7}{15},\frac{1}{2})$}\\
0&\text{otherwise.}
\end{cases}
\]
Thus $m(x)-\tilde{m}(x)= 2$ in a neighborhood of $x= 0$, and $a=1$ does not satisfy the hypothesis of Proposition \ref{constructH} or Theorem \ref{kathy}.  However, the values $a=2$ and $a=3$ satisfy the hypothesis of Proposition \ref{constructH}, and we can construct filters corresponding to these two values.

For $a=2$, we can take 
\[
H=[h_{i,j}]=\sqrt{2}\begin{pmatrix}
\chi_{[-\frac{2}{15},\frac{2}{15})\cup \pm [\frac{1}{4},\frac{4}{15})}&0&0\\
0&\chi_{[-\frac{1}{15},\frac{1}{15})}&0\\
\chi_{\pm (\frac{7}{15},\frac{1}{2}]}&0&0
\end{pmatrix},
\]
which is low-pass of rank $a=2$ for dilation by $2$.  A different low-pass filter of rank $a=2$ for this $m$ is 
\[
H=[h_{i,j}]=\sqrt{2}\begin{pmatrix}
\chi_{[-\frac{2}{15},\frac{2}{15})\cup \pm [\frac{1}{4},\frac{4}{15})}&0&0\\
0&\chi_{[-\frac{1}{30},\frac{1}{30})}&\chi_{\pm [\frac{1}{30},\frac{1}{15})}\\
\chi_{\pm (\frac{13}{30},\frac{1}{2}]}&0&0
\end{pmatrix}.
\]

The following family satisfies the low-pass condition of rank $a=3$ for dilation by $2$, with the same choice of $m$:
\[
H=[h_{i,j}]=\sqrt{2}\begin{pmatrix}
\chi_{[-\frac{2}{15},\frac{2}{15})\cup \pm [\frac{1}{4},\frac{4}{15})}&0&0\\
0&\chi_{[-\frac{1}{15},\frac{1}{15})}&0\\
0&0&\chi_{ [-\frac{1}{30},\frac{1}{30})}
\end{pmatrix}.
\]

In all of these cases, the infinite product $\prod_{j=1}^{\infty}[2^{-1/2}H(2^{-j}x)]$ gives functions whose inverse Fourier transforms in $L^2(\mathbb R)$ generate a shift invariant subspace with multiplicity function a degenerate form of the original $m$. \end{example}

 \end{document}